\documentclass[12pt,reqno]{amsart}
\usepackage{amssymb,verbatim,enumerate,ifthen}
\usepackage[mathscr]{eucal}
\usepackage[utf8]{inputenc}
\usepackage[T1]{fontenc}
\allowdisplaybreaks[1]
\def\N{\mathbb{N}}
\def\R{\mathbb{R}}

\def\Cdot{\!\cdot\!}

\newtheorem{theorem}{Theorem}
\newtheorem*{theorem*}{Theorem}
\def\Thm#1#2{\ifthenelse{\equal{#1}{*}}{\begin{theorem*}#2\end{theorem*}}
             {\begin{theorem}\label{T#1}#2\end{theorem}}}
\newtheorem{Atheorem}{Theorem}

\def\thm#1{Theorem~\ref{T#1}}

\newtheorem{proposition}[theorem]{Proposition}
\newtheorem*{proposition*}{Proposition}
\def\Prp#1#2{\ifthenelse{\equal{#1}{*}}{\begin{proposition*}#2\end{proposition*}}
             {\begin{proposition}\label{P#1}#2\end{proposition}}}

\newtheorem{corollary}[theorem]{Corollary}
\newtheorem*{corollary*}{Corollary}
\def\Cor#1#2{\ifthenelse{\equal{#1}{*}}{\begin{corollary*}#2\end{corollary*}}
             {\begin{corollary}\label{C#1}#2\end{corollary}}}

\newtheorem{lemma}[theorem]{Lemma}
\newtheorem*{lemma*}{Lemma}
\def\Lem#1#2{\ifthenelse{\equal{#1}{*}}{\begin{lemma*}#2\end{lemma*}}
             {\begin{lemma}\label{L#1}#2\end{lemma}}}
\def\lem#1{Lemma~\ref{L#1}}
\newtheorem{Alemma}{Lemma}

\theoremstyle{definition}
\newtheorem{remark}[theorem]{Remark}
\newtheorem*{remark*}{Remark}
\def\Rem#1#2{\ifthenelse{\equal{#1}{*}}{\begin{remark*}\rm #2\end{remark*}}
             {\begin{remark}\label{R#1}\rm #2\end{remark}}}

\newtheorem{example}[theorem]{Example}
\newtheorem*{example*}{Example}
\def\Exa#1#2{\ifthenelse{\equal{#1}{*}}{\begin{example*}\rm #2\end{example*}}
             {\begin{example}\label{Ex#1}\rm #2\end{example}}}

\def\eq#1{{\rm(\ref{E#1})}}
\def\Eq#1#2{\ifthenelse{\equal{#1}{*}}
  {\begin{equation*}\begin{aligned}#2\end{aligned}\end{equation*}}
  {\begin{equation}\begin{aligned}\label{E#1}#2\end{aligned}\end{equation}}}

\def\diag{\mathop{\hbox{\rm diag}}\nolimits}
\def\conv{\mathop{\hbox{\rm conv}}\nolimits}
\def\id{\mathop{\hbox{\rm id}}\nolimits}

\begin{document}
%\begin{flushright}
%\textit{Submitted to:} 
%\end{flushright}
\vspace{5mm}

\date{\today}

\title[Equality of generalized Bajraktarević means]{On the equality problem of generalized Bajraktarević means}

\author[R.\ Gr\"unwald]{Rich\'ard Gr\"unwald}
\author[Zs. P\'ales]{Zsolt P\'ales}
\address{Institute of Mathematics, University of Debrecen, 
H-4002 Debrecen, Pf. 400, Hungary}
\email{richard.grunwald96@gmail.com, pales@science.unideb.hu}

\thanks{The research of the first author was supported by the \'UNKP-18-2 New National Excellence Program of 
	the Ministry of Human Capacities. The research of the second author was supported by the Hungarian Scientific Research Fund (OTKA) Grant
	K-111651 and by the EFOP-3.6.1-16-2016-00022 project. This project is co-financed by the European Union and the European Social Fund.}

\dedicatory{Dedicated to the 95th birthday of Professor János Aczél}

\begin{abstract}
The purpose of this paper is to investigate the equality problem of generalized Bajraktarević means, i.e., to solve the functional equation
\begin{equation}\label{E0}\tag{*}
	f^{(-1)}\bigg(\frac{p_1(x_1)f(x_1)+\dots+p_n(x_n)f(x_n)}{p_1(x_1)+\dots+p_n(x_n)}\bigg)=g^{(-1)}\bigg(\frac{q_1(x_1)g(x_1)+\dots+q_n(x_n)g(x_n)}{q_1(x_1)+\dots+q_n(x_n)}\bigg),
\end{equation}
which holds for all $x=(x_1,\dots,x_n)\in I^n$, where $n\geq 2$, $I$ is a nonempty open real interval, the unknown functions $f,g:I\to\R$ are strictly monotone, $f^{(-1)}$ and $g^{(-1)}$ denote their generalized left inverses, respectively, and $p=(p_1,\dots,p_n):I\to\R_{+}^n$ and $q=(q_1,\dots,q_n):I\to\R_{+}^n$ are also unknown functions. This equality problem in the symmetric two-variable (i.e., when $n=2$) case was already investigated and solved under sixth-order regularity assumptions by Losonczi in 1999. In the nonsymmetric two-variable case, assuming three times differentiability of $f$, $g$ and the existence of $i\in\{1,2\}$ such that either $p_i$ is twice continuously differentiable and $p_{3-i}$ is continuous on $I$, or $p_i$ is twice differentiable and $p_{3-i}$ is once differentiable on $I$, we prove that \eq{0} holds if and only if there exist four constants $a,b,c,d\in\R$ with $ad\neq bc$ such that 
\Eq{*}{
        cf+d>0,\qquad 
	g=\frac{af+b}{cf+d},\qquad\mbox{and}\qquad q_\ell=(cf+d)p_\ell\qquad (\ell\in\{1,\dots,n\}).
}
In the case $n\geq 3$, we obtain the same conclusion with weaker regularity assumptions. Namely, we suppose that $f$ and $g$ are three times differentiable, $p$ is continuous and there exist $i,j,k\in\{1,\dots,n\}$ with $i\neq j\neq k\neq i$ such that $p_i,p_j,p_k$ are differentiable. 
\end{abstract}

\maketitle

\section{Introduction}

Throughout this paper, the symbols $\R$ and $\R_+$ will stand for the sets of real and positive real numbers, respectively, and $I$ will always denote a nonempty open real interval. In theory of quasi-arithmetic means the characterization of the equality of means with different generators is a basic problem which was completely solved in the book \cite{HarLitPol34}. Using this characterization, the homogeneous quasi-arithmetic means can also be found: they are exactly the power means and the geometric mean. In \cite{Baj58} (cf.\ also \cite{Baj63}) Bajraktarević introduced a new generalization of quasi-arithmetic means by adding a weight function to the formula of quasi-arithmetic means. He also described the equality of such means (called Bajraktarević means since then) in the at least $3$-variable setting assuming three times differentiability. Daróczy and Losonczi in \cite{DarLos70}, later Daróczy and Páles in \cite{DarPal82} arrived at the same conlusion with first-order differentiability and without differentiability, respectively, but assuming the equality for all $n\in\N$. As an application of the characterization of the equality, Aczél and Daróczy in \cite{AczDar63c} determined the homogeneous Bajraktarević means that include Gini means which were introduced by Gini in \cite{Gin38}. Losonczi in \cite{Los99} described the equality of two-variable Bajraktarević means under sixth-order regularity assumptions and an algebraic condition which was later removed in \cite{Los06b}. Using these results, the homogeneous two-variable means were also determined by Losonczi \cite{Los07a}, \cite{Los07b}.

The purpose of this paper is to extend the definition of Bajraktarević means in a nonsymmetric way by replacing each appearance of the weight function by a possibly different one. We also take strictly monotone functions instead of strictly monotone and continuous ones.

Given a subset $S\subseteq\R$, the smallest convex set containing $S$, which is identical to the smallest interval containing $S$, will be denoted by $\conv(S)$.
For our definition of generalized Bajraktarević means, we shall need the following lemma about the existence and properties of the left inverse of strictly monotone (but not necessarily continuous) functions.

\Lem{SMF}{Let $f:I\to\R$ be a strictly monotone function. Then there exists a uniquely determined monotone function $g:\conv(f(I))\to I$ such that $g$ is the left inverse of $f$, i.e.,
\Eq{SMF1}{
   (g\circ f)(x)=x \qquad(x\in I).
}
Furthermore, $g$ is monotone in the same sense as $f$, continuous,
\Eq{SMF2}{
  (f\circ g)(y)=y \qquad(y\in f(I)),
}
and
\Eq{SMF3}{
  \liminf_{x\to g(y)}f(x)\leq y \leq \limsup_{x\to g(y)}f(x) 
  \qquad(y\in\conv(f(I))).
}
Thus, if $f$ is lower (resp.\ upper) semicontinuous at $g(y)$, then $f\circ g(y)\leq y$ (resp.\ $y\leq f\circ g(y)$).}

\begin{proof} Without loss of generality, we may assume that $f:I\to\R$ is a strictly increasing function. Then $f:I\to f(I)$ is a bijection. The interval $I$ is open, therefore, $f$ has a left and a right limit at every point $x\in I$, which will be denoted by $f_-(x)$ and $f_+(x)$, respectively. We introduce the notation $J_x:=[f_-(x),f_+(x)]$, where $x\in I$. Then, for all elements $u<x<v\in I$, we have that
\Eq{*}{
   f_+(u)<f_-(x)\leq f(x)\leq f_+(x)<f_-(v).
}
From these inequalities, it follows that $f(x)\in J_x$ holds for all $x\in I$ and $J_x\cap J_u=\emptyset$ whenever $u$ is distinct from $x$. 

The convex hull of $f(I)$ is the smallest interval $J\subseteq\R$ containing $f(I)$. The opennes of $I$ implies that $\inf f(I),\sup f(I)\not\in J$, hence $J:=\,]\inf f(I),\sup f(I)[\,$. We show that 
\Eq{JJ}{
   J=\bigcup_{x\in I} J_x.
}
If $x\in I$, then, for all $u<x$, we have $f_-(x)>f_+(u)=\inf_{u<t} f(t)\geq \inf f(I)$.
Similarly, $f_+(x)<\sup f(I)$, therefore, $J_x\subseteq J$. This proves the inclusion $\supseteq$ in \eq{JJ}. To prove the reversed inclusion in \eq{JJ}, let $y\in J$. Define
\Eq{*}{
  x:=\sup\{u\in I\mid f(u)\leq y\}.
}
Then, for all $n\in\N$, there exists $u_n\in I$ such that $x-\frac1n<u_n$ and $f(u_n)\leq y$. Thus, $u_n\leq x$ and hence $u_n$ tends to $x$ as $n\to\infty$. Therefore,
\Eq{*}{
   f_-(x)\leq\limsup_{n\to\infty} f(u_n)\leq y.
}
On the other hand, let $u_n\in I$ be an arbitrary sequence converging to $x$ such that $x<u_n$. Then $y<f(u_n)$, whence we obtain
\Eq{*}{
	y\leq\lim_{n\to\infty} f(u_n)=f_+(x).
} 
The above inequalities imply that $y\in J_x$, which completes the proof of the inclusion $\subseteq$ in \eq{JJ}.

Let $y\in J=\conv(f(I))$ be an arbitrarily fixed element. Then there exists a uniquely determined element $x\in I$ such that $y\in J_x$, hence we define the function $g:J\to I$ by the prescription $g(y):=x$. 

Therefore, if $x\in I$ is an arbitrary element, then it is obvious that $f(x)\in J_x$ and hence $g(f(x))=x$. Thus, equation \eq{SMF1} is valid for all $x\in I$. 

To see that $g$ is nondecreasing, let $y_1<y_2$ be arbitrary elements of $J$. Then there exist elements $x_1,x_2\in I$ such that $y_i\in J_{x_i}$. If $x_2$ were strictly smaller than $x_1$, then we would have
\Eq{*}{
  y_2\leq f_+(x_2)<f_-(x_1)\leq y_1.
}
This contradiction shows that $g(y_1)=x_1\leq x_2=g(y_2)$.

To prove that $g$ is continuous, let $y\in J$ and choose $\varepsilon>0$ so that $g(y)\pm\varepsilon$ be in $I$. Define $W_\varepsilon:=\,]f_-(g(y)-\varepsilon),f_+(g(y)+\varepsilon)[\,$. Then
\Eq{*}{
   f_-(g(y)-\varepsilon)<f_-(g(y))\leq y\leq f_+(g(y))<f_+(g(y)+\varepsilon),
}   
hence $W_\varepsilon$ is neighborhood of $y$. 
By the monotonicity of $g$, for $w\in W_\varepsilon$, we have that
\Eq{*}{
  g(y)-\varepsilon=g\big(f_-(g(y)-\varepsilon)\big)
  \leq g(w) \leq g\big(f_+(g(y)+\varepsilon)\big)
  =g(y)+\varepsilon,
}
which yields that $g$ is continuous at $y$.

If $y\in f(I)$, then there exists a uniquely determined element $x\in I$ such that $f(x)=y$ and hence, using \eq{SMF1}, we get that 
\Eq{*}{
	(f\circ g)(y)=f((g\circ f)(x))=f(x)=y,
}
which shows that \eq{SMF2} holds for all $y\in f(I)$.

To see that \eq{SMF3} is valid, let $y\in J$. By the definition of $g(y)$, there exists a unique element $v\in I$ such that $y\in J_v$ and $g(y)=v$. Then, for all $x<v=g(y)$, we have 
\Eq{*}{
   f(x)\leq f_+(x)<f_-(v)\leq y.
}
Therefore, upon taking the left limit $x\to v-0$, we get
\Eq{*}{
   \liminf_{x\to g(y)}f(x)=\lim_{x\to g(y)-0}f(x)\leq y,
}
which proves the left hand side inequality in \eq{SMF3}. The verification of the right hand side inequality is completely analogous, therefore it is omitted.

Finally, we prove the uniqueness of $g$. Assume that $h:J\to I$ is a nondecreasing function which is the left inverse of $f$. We are going to show that $h$ coincides with $g$ on $J$. Let $y\in J$ be arbitrary. Then there exists $x\in I$ such that $f_-(x)\leq y\leq f_+(x)$ and $g(y)=x$. Let $(x_n)$ be a strictly increasing and $(x_n')$ be a strictly decreasing sequence converging to $x$. Then, for all $n\in\N$, we have
\Eq{*}{
  f(x_n)<f_-(x)\leq y\leq f_+(x)<f(x_n').
}
By the monotonicity of $h$, it follows that
\Eq{*}{
  x_n=(h\circ f)(x_n)\leq h(y)\leq (h\circ f)(x_n')=x_n'.
}
Taking the limit $n\to\infty$, we arrive at
\Eq{*}{
   x\leq h(y)\leq x,
}
which proves that $h(y)=x=g(y)$.
\end{proof}

The function $g$ described in the above lemma is called the \emph{generalized left inverse of the strictly monotone function $f:I\to\R$} and is denoted by $f^{(-1)}$. It is clear from \eq{SMF1} and \eq{SMF2} that the restriction of $f^{(-1)}$ to $f(I)$ is the inverse of $f$ in the standard sense. Therefore, $f^{(-1)}$ is the continuous and monotone extension of the inverse of $f$ to the smallest interval containing the range of $f$.

Given a strictly monotone function $f:I\to\R$ and an $n$-tuple of positive valued functions $p=(p_1,\dots,p_n):I\to\R_{+}^n$, we introduce the \emph{$n$-variable generalized Bajraktarević mean} $A_{f,p}:I^n\to I$ by the following formula:
\Eq{BM}{
	A_{f,p}(x):=f^{(-1)}\bigg(\frac{p_1(x_1)f(x_1)+\dots+p_n(x_n)f(x_n)}{p_1(x_1)+\dots+p_n(x_n)}\bigg) \qquad (x=(x_1,\dots,x_n)\in I^n),
}
and, to simplify the notations, we will use the following definition:
\Eq{Rfp}{
	R_{f,p}(x)
	:=\frac{p_1(x_1)f(x_1)+\dots+p_n(x_n)f(x_n)}{p_1(x_1)+\dots+p_n(x_n)}.
}

\Thm{BM1}{Let $f:I\to\R$ be strictly monotone and $p=(p_1,\dots,p_n):I\to\R_{+}^n$. Then the function $A_{f,p}:I^n\to I$ given by \eq{BM} is well-defined and it is a mean, that is, 
\Eq{MV}{
	\min(x)\leq A_{f,p}(x)\leq \max(x)\qquad (x=(x_1,\dots,x_n)\in I^n).
}
}

\begin{proof}
We may assume that $f$ is strictly increasing (in the decreasing case the proof is completely similar). To show that, for all $x=(x_1,\dots,x_n)\in I^n$, the formula for $A_{f,p}(x)$ is well-defined and \eq{MV} holds, consider the ratio $R_{f,p}(x)$.

Due to the positivity of the values of $p_i(x_i)$, we can see that $R_{f,p}(x)$ is a convex combination of the values $f(x_1),\dots,f(x_n)$, therefore,
\Eq{MV1}{
  f(\min(x))=\min(f(x_1),\dots,f(x_n))
  \leq R_{f,p}(x)\leq \max(f(x_1),\dots,f(x_n))=f(\max(x)).
}
This shows that $R_{f,p}(x)$ is an element of $\conv(f(I))$, which is the domain of $f^{(-1)}$ and hence $A_{f,p}(x)=f^{(-1)}(R_{f,p}(x))$ is well-defined. Furthermore, using that $f^{(-1)}$ is nondecreasing and is the left inverse of $f$, the inequalities in \eq{MV1} yield
\Eq{*}{
  \min(x)=f^{(-1)}(f(\min(x)))\leq f^{(-1)}(R_{f,p}(x))\leq f^{(-1)}(f(\max(x)))=\max(x).
}
This finally proves the mean value inequalities stated in \eq{MV}.
\end{proof}

\Thm{BM2}{Let $f:I\to\R$ be strictly increasing and $p=(p_1,\dots,p_n):I\to\R_{+}^n$. Then, for all $x=(x_1,\dots,x_n)\in I^n$, the equality $y=A_{f,p}(x)$ holds if and only if 
	\Eq{MV2}{
	\sum_{i=1}^{n}p_i(x_i)(f(z)-f(x_i))
	\begin{cases}
	 <0 & \mbox{ for } z\in I,\, z<y,\\[2mm]
	 >0 & \mbox{ for } z\in I,\, z>y.
	\end{cases}
	}
If $f$ is strictly decreasing, then the inequalities \eq{MV2} hold with reversed inequality sign.
}

\begin{proof} Assume that $f:I\to\R$ is strictly increasing, let $x=(x_1,\dots,x_n)\in I^n$ and $y:=A_{f,p}(x)$. If $z<y$, then  $f(z)<R_{f,p}(x)$, because in the opposite case we would have $f(z)\geq R_{f,p}(x)$ which implies $z=f^{(-1)}(f(z))\geq f^{(-1)}(R_{f,p}(x))=A_{f,p}(x)=y$, contradicting the choice of $z$. Rearranging the inequality $f(z)<R_{f,p}(x)$, it easily follows that \Eq{*}{
\sum_{i=1}^{n}p_i(x_i)(f(z)-f(x_i))<0.
}
In the case $z>y$, we get $f(z)>R_{f,p}(x)$, which implies the second inequality in \eq{MV2}.

Observe that the function 
\Eq{*}{
  z\mapsto \varphi(z):=\sum_{i=1}^{n}p_i(x_i)(f(z)-f(x_i))
}
is strictly increasing. Therefore, it changes sign at at most one point in $I$.
If \eq{MV2} holds for $y$, then $\varphi$ changes sign at $y$. On the other hand, as we have seen it above, $\varphi$ also changes sign at $A_{f,p}(x)$. Hence $y=A_{f,p}(x)$ must hold.
\end{proof}

\Cor{BM2}{Let $f:I\to\R$ be continuous, strictly monotone, and $p=(p_1,\dots,p_n):I\to\R_{+}^n$. Then, for all $x=(x_1,\dots,x_n)\in I^n$, the value $y=A_{f,p}(x)$ is the unique solution of the equation
	\Eq{MV3}{
	\sum_{i=1}^{n}p_i(x_i)(f(y)-f(x_i))=0.
	}
}

\begin{proof}
	The function 
	\Eq{*}{
		y\mapsto \varphi(y):=\sum_{i=1}^{n}p_i(x_i)(f(y)-f(x_i))
	}
	is strictly monotone and continuous. Therefore, it vanishes at most one point in $I$. Applying \thm{BM2}, we obtain that $\varphi$ changes sign at $y=A_{f,p}(x)$. Thus, using that $\varphi$ is continuous, $\varphi$ vanishes at $y=A_{f,p}(x)$.
\end{proof}

The next result establishes a sufficient condition for the equality of the $n$-variable generalized Bajraktarević means. We will call this situation the canonical case of the equality.

\Thm{BM3}{
	Let $f,g:I\to\R$ be strictly monotone and $p=(p_1,\dots,p_n):I\to\R_{+}^n,\,q=(q_1,\dots,q_n):I\to\R_{+}^n$. If there exist $a,b,c,d\in\R$ with $ad\neq bc$ such that 
\Eq{fg}{
        cf+d>0,\qquad
	g=\frac{af+b}{cf+d},\qquad\mbox{and}\qquad q_i=(cf+d)p_i\qquad (i\in\{1,\dots,n\})
}
hold on $I$, then the $n$-variable generalized Bajraktarević means $A_{f,p}$ and $A_{g,q}$ are identical on $I^n$.
}

\begin{proof}
	Let $x=(x_1,\dots,x_n)\in I^n$ be arbitrary. Using the formulas \eq{fg}, we obtain that
	\Eq{*}{
	\sum_{i=1}^{n}q_i(x_i)(g(z)&-g(x_i))\\&=\sum_{i=1}^{n}(cf(x_i)+d)p_i(x_i)\left(\frac{(af(z)+b)(cf(x_i)+d)-(af(x_i)+b)(cf(z)+d)}{(cf(x_i)+d)(cf(z)+d)}\right)\\&=\frac{ad-bc}{cf(z)+d}\Big(\sum_{i=1}^{n}p_i(x_i)(f(z)-f(x_i))\Big).
}
It shows that $\sum_{i=1}^{n}q_i(x_i)(g(z)-g(x_i))$ changes sign at $y$ if and only if $\sum_{i=1}^{n}p_i(x_i)(f(z)-f(x_i))$ changes sign at $y$. Hence, applying \thm{BM2}, $A_{f,p}(x)=A_{g,q}(x)$ holds. The element $x$ being arbitrary in $I^n$, we get the statement of the theorem.
\end{proof}

With the aid of the following lemma, we can reduce the regularity assumptions in our statements. For the formulation of this and the subsequent results, 
we define the \emph{diagonal} $\diag(I^n)$ of $I^n$ and the map $\Delta_n:I\to\diag(I^n)$ by 
\Eq{*}{
\diag(I^n):=\{(x,\dots,x)\in\R^n\mid x\in I\}\qquad\mbox{and}\qquad
  \Delta_n(x):=(x,\dots,x) \qquad(x\in I).
}
For all $i\in\{1,\dots,n\}$, let $e_i\in\R^n$ denote the $i$th vector of the standard base of $\R^n$, i.e., let $e_i:=(\delta_{ij})_{j=1}^n$, where $\delta$ stands for the Kronecker symbol.

Given $p=(p_1,\dots,p_n):I\to\R_{+}^n$ and $q=(q_1,\dots,q_n):I\to\R_{+}^n$, we will also use the following notations:
\Eq{*}{
	p_0:=p_1+\dots+p_n,\qquad q_0:=q_1+\dots+q_n,\qquad\mbox{and}\qquad r_0:=\frac{q_0}{p_0}.
} 

\Lem{reg}{
        Let $f,g:I\to\R$ be continuous strictly monotone functions, $n\geq2$, and $p=(p_1,\dots,p_n):I\to\R_{+}^n,\,q=(q_1,\dots,q_n):I\to\R_{+}^n$. Assume that there exists an open set $U\subseteq I^n$ containing $\diag(I^n)$ such that $A_{f,p}=A_{g,q}$ holds on $U$.  Then the following two assertions hold.
        \begin{enumerate}[(i)]
         \item For all $i\in\{1,\dots,n\}$, the function $p_i$ is continuous on $I$ if and only if the function $q_i$ is continuous on $I$.
         \item Let $k\in\N$. Assume that $f,g:I\to\R$ are $k$ times differentiable (resp.\ $k$ times continuously differentiable) functions on $I$ with nonvanishing first derivatives. Then, for all $i\in\{1,\dots,n\}$, the function $p_i$ is $k$ times differentiable (resp.\ $k$ times  continuously differentiable) on $I$ if and only if $q_i$ is $k$ times differentiable (resp.\ $k$ times continuously differentiable) on $I$.
        \end{enumerate}
}

\begin{proof} In what follows, we will prove that the regularity properties possessed by $p_i$ are transferred to the corresponding $q_i$. The reversed statements can similarly be verified.

For $i\in\{1,\dots,n\}$, denote 
\Eq{*}{
  U_i:=\{(x,y)\in I^2\mid \Delta_n(x)+(y-x)e_i\in U\}.
}
Then $U_i$ is an open set containing $\diag(I^2)$. By our assumption, we have that, for all $(x,y)\in U_i$,
\Eq{*}{
  A_{g,q}(\Delta_n(x)+(y-x)e_i)=A_{f,p}(\Delta_n(x)+(y-x)e_i).
}
This is equivalent to the following equality
\Eq{ee}{
   \frac{(q_0(x)-q_i(x))g(x)+q_i(y)g(y)}{q_0(x)-q_i(x)+q_i(y)}
   =\big(g\circ f^{-1}\big)\bigg(\frac{(p_0(x)-p_i(x))f(x)+p_i(y)f(y)}{p_0(x)-p_i(x)+p_i(y)}\bigg) \qquad((x,y)\in U_i).
}
Observe that, for $x,y\in I$ with $x\neq y$, the inequalities $p_i(x)<p_0(x)$ and $f(x)\neq f(y)$ imply that 
\Eq{*}{\frac{
		(p_0(x)-p_i(x))f(x)+p_i(y)f(y)}{p_0(x)-p_i(x)+p_i(y)}\neq f(y).
}
Therefore, 
\Eq{*}{
	\big(g\circ f^{-1}\big)\bigg(\frac{(p_0(x)-p_i(x))f(x)+p_i(y)f(y)}{p_0(x)-p_i(x)+p_i(y)}\bigg)\neq g(y).
}
Thus, solving equation \eq{ee} with respect to $q_i(y)$, we get
\Eq{qi}{
  q_i(y)=(q_0(x)-q_i(x))\frac{(g\circ f^{-1})\big(\frac{(p_0(x)-p_i(x))f(x)+p_i(y)f(y)}{p_0(x)-p_i(x)+p_i(y)}\big)-g(x)}{g(y)-(g\circ f^{-1})\big(\frac{(p_0(x)-p_i(x))f(x)+p_i(y)f(y)}{p_0(x)-p_i(x)+p_i(y)}\big)}
  \qquad((x,y)\in U_i,\,x\neq y).
}
Let $x_0\in I$ be an arbitrarily fixed point. The pair $(x_0,x_0)$ is an interior point of $U_i$, therefore, there exists $x\in I\setminus\{x_0\}$ such that $(x,x_0)\in U_i$. Then the set 
\Eq{*}{
  V_i:=\{y\in I\mid (x,y)\in U_i,\,x\neq y\} 
}
is a neighborhood of $x_0$ on which we have the equality \eq{qi} for $q_i$. 

Provided that $f$ and $g$ are continuous on $I$ and $p_i$ is continuous at $x_0$, it follows that $g\circ f^{-1}$ is continuous on $f(I)$ and hence the mapping 
\Eq{map}{
  y\mapsto \big(g\circ f^{-1}\big)\bigg(\frac{(p_0(x)-p_i(x))f(x)+p_i(y)f(y)}{p_0(x)-p_i(x)+p_i(y)}\bigg)
}
is continuous at $x_0$. This shows that the right hand side of \eq{qi} is a continuous function of $y$ at $x_0$ and hence $q_i$ is continuous at $x_0$.
This proves the first assertion.

Provided that, for some $k\in\N$, the functions $f,g:I\to\R$ are $k$ times differentiable (resp.\ $k$ times continuously differentiable) on $I$ with nonvanishing first derivatives and that $p_i$ is $k$ times differentiable (resp.\ $k$ times continuously differentiable) at $x_0$, it follows, by the standard calculus rules, that $g\circ f^{-1}$ is $k$ times differentiable (resp.\ $k$ times continuously differentiable) and hence the mapping \eq{map} is also $k$ times differentiable (resp.\ $k$ times continuously differentiable) at $x_0$. This implies that the right hand side of \eq{qi} is a $k$ times differentiable (resp.\ $k$ times continuously differentiable) function of $y$ at $x_0$ and hence $q_i$ is $k$ times differentiable (resp.\ $k$ times continuously differentiable) at $x_0$. This proves the second statement.
\end{proof}

The following theorem is of basic importance for our investigations. 

\Thm{BM4}{
	Let $f,g:I\to\R$ be continuous, strictly monotone and $p=(p_1,\dots,p_n):I\to\R_{+}^n$ be continuous function on $I$. Let further $q=(q_1,\dots,q_n):I\to\R_{+}^n$. Assume that there exists an open set $U\subseteq I^n$ containing the $\diag(I^n)$ such that $A_{f,p}=A_{g,q}$ holds on $U$ and that there exist $a,b,c,d\in\R$ with $ad\neq bc$ and a nonempty open subinterval $J$ of $I$ such that \eq{fg} holds on $J$. Then $q$ is continuous on $I$ and \eq{fg} is also valid on $I$.
}

\begin{proof}
First of all, using \lem{reg} and the continuity of $f$, $g$ and $p$, it is clear that $q$ is continuous on $I$.

Assume that $A_{f,p}=A_{g,q}$ holds on some open set $U$ containing the $\diag(I^n)$ and for some constants $a,b,c,d\in\R$ with $ad\neq bc$ there exists a nonempty open subinterval $J$ of $I$ such that \eq{fg} holds on $J$. We may assume that $J$ is a maximal subinterval of $I$ with this property. To complete the proof, we have to show that $J=I$. To the contrary, suppose that $J\neq I$. Then one of the strict inequalities 
\Eq{inf}{
	\inf I<\inf J=:\alpha \qquad\mbox{or}\qquad\sup J<\sup I
}
must be valid. We may suppose that first inequality in \eq{inf} holds. Hence, due to the continuity of $f$, $p_1$, and $q_1$ at $\alpha$, it follows from \eq{fg} that $q_1(\alpha)=(cf(\alpha)+d)p_1(\alpha)$. Therefore, $q_1(\alpha)>0$ implies that $cf(\alpha)+d>0$. Consequently, using the continuity of all functions, for all $x\in\bar{J}:=J\cup\{\alpha\}$, we get that
\Eq{*}{
	cf(x)+d>0,\qquad
	g(x)=\frac{af(x)+b}{cf(x)+d}\qquad\mbox{and}\qquad 
	q_i(x)=(cf(x)+d)p_i(x) \qquad (i\in\{1,\dots,n\})
}
are valid. By the continuity of $f$, there is an element $\bar\alpha\in I$ with
$\bar\alpha<\alpha$ such that $cf(x)+d>0$ for all $x\in\bar{I}:=\,]\bar\alpha,\alpha]\cup J$. Define the functions $\bar{g}:\bar{I}\to\R$ and $\bar{q}:\bar{I}\to\R_+^n$ by 
\Eq{gq}{
	\bar{g}(x):=\frac{af(x)+b}{cf(x)+d} \qquad\mbox{and}\qquad
	\bar{q_i}(x):=(cf(x)+d)p_i(x) \qquad(x\in\bar{I},\,i\in\{1,\dots,n\}).
}
Thus, for all $x\in\bar{J}$, the equations 
\Eq{gg}{
  g(x)=\bar{g}(x)\qquad\mbox{and}\qquad q_i(x)=\bar{q_i}(x) 
  \qquad (i\in\{1,\dots,n\})
}
hold. On the other hand, the maximality property of $J$ implies that there is no $\beta<\alpha$ such that \eq{gg} is valid for all $x\in\,]\beta,\alpha]\cup J$. Furthermore, the equality $A_{f,p}=A_{g,q}$ on $U$ and \thm{BM3} applied to the conditions \eq{gq} yield that 
\Eq{AAA}{
  A_{g,q}(x)=A_{f,p}(x)=A_{\bar{g},\bar{q}}(x)
} 
is also valid for all $x\in \bar{U}:=(\bar{I})^n\cap U$. The point $(\alpha,\dots,\alpha)$ is an interior point of $\bar{U}$, therefore, there exists $r>0$ such that $\,]\alpha-r,\alpha+r[^n\,\subseteq\bar{U}$ and hence \eq{AAA} holds for all $x\in\,]\alpha-r,\alpha+r[^n$.

In what follows, we assume that $g$ is strictly increasing and hence $\bar{g}$ must be also strictly increasing. The functions $g$ and $\bar{g}$ are identical 
on $[\alpha,\alpha+r[$, therefore, their inverses are also equal on $[g(\alpha),g(\alpha+r)[$.

The following claim will be useful for the rest of the proof.

\noindent\textbf{Claim.} \textit{If $x=(x_1,\dots,x_n)\in\,]\alpha-r,\alpha+r[^n$ such that $\alpha\leq A_{g,q}(x)$, then}
\Eq{Claim}{
  \frac{q_1(x_1)g(x_1)+\dots+q_n(x_n)g(x_n)}{q_1(x_1)+\dots+q_n(x_n)}
  =\frac{\bar{q}_1(x_1)\bar{g}(x_1)+\dots+\bar{q}_n(x_n)\bar{g}(x_n)}{\bar{q}_1(x_1)+\dots+\bar{q}_n(x_n)}.
}
Indeed, the condition on $x$ implies that $\alpha \leq A_{g,q}(x)\leq\max(x)<\alpha+r$ also holds, hence $g(A_{g,q}(x))=\bar{g}(A_{g,q}(x))$. On the other hand, in view of \eq{AAA}, we have the equality $A_{g,q}(x)=A_{\bar{g},\bar{q}}(x)$. Therefore, $g(A_{g,q}(x))=\bar{g}(A_{\bar{g},\bar{q}}(x))$, which implies the equation \eq{Claim}.

Let $y_0\in\,]\alpha,\alpha+r[$ be fixed. Then the inequality $g(\alpha)<g(y_0)$ implies that
\Eq{*}{
  g(\alpha)<\frac{q_i(\alpha)g(\alpha)+(q_0-q_i)(y_0)g(y_0)}{q_i(\alpha)+(q_0-q_i)(y_0)} \qquad(i\in\{1,\dots,n\}).
}
Now, by the continuity of the functions $g,q_1,\dots,q_n$, we can find a positive number $\delta_0:=\delta(y_0)<\min(y_0-\alpha,\alpha+r-y_0)<r$ such that, for all $x\in\,]\alpha-\delta_0,\alpha]$ and $y\in\,]y_0-\delta_0,y_0+\delta_0[$, 
\Eq{gal}{
  g(\alpha)\leq\frac{q_i(x)g(x)+(q_0-q_i)(y)g(y)}{q_i(x)+(q_0-q_i)(y)} 
  \qquad(i\in\{1,\dots,n\}).
}
Applying the inverse of $g$ side by side to this inequality, it follows that
$\alpha\leq A_{g,q}(x_1,\dots,x_n)$, where $x_i:=x$ and $x_j:=y$ for all $j\in\{1,\dots,n\}\setminus\{i\}$. Therefore, in view of the Claim above and the equality \eq{gg}, for all $x\in\,]\alpha-\delta_0,\alpha]$ and $y\in\,]y_0-\delta_0,y_0+\delta_0[$, we have that
\Eq{*}{
  \frac{q_i(x)g(x)+(q_0-q_i)(y)g(y)}{q_i(x)+(q_0-q_i)(y)}
  =\frac{\bar{q}_i(x)\bar{g}(x)+(q_0-q_i)(y)g(y)}{\bar{q}_i(x)+(q_0-q_i)(y)}
  \qquad(i\in\{1,\dots,n\}).
}
This equality can be rewritten as
\Eq{3tag}{
  q_i(x)\bar{q}_i(x)(g(x)\!-\!\bar{g}(x))
  +(q_0\!-\!q_i)(y)(q_i(x)g(x)\!-\!\bar{q}_i(x)\bar{g}(x))
  +(q_0\!-\!q_i)(y)g(y)(\bar{q}_i(x)\!-\!q_i(x))=0.
}

Consider the sets
\Eq{*}{
  S:=\{x\in\,]\alpha-r,\alpha[\,\colon g(x)\neq\bar{g}(x)\},\quad
  S_i:=\{x\in\,]\alpha-r,\alpha[\,\colon q_i(x)\neq\bar{q}_i(x)\},
  \quad(i\in\{1,\dots,n\}).
}
In the next step we show that 
\Eq{SS}{
   S\,\cap \,]\alpha-\delta_0,\alpha[\,
   = S_i\,\cap \,]\alpha-\delta_0,\alpha[\, \qquad(i\in\{1,\dots,n\}).
}
If $x\in\,]\alpha-\delta_0,\alpha[\,\setminus S$, then $g(x)=\bar{g}(x)$. Using this, \eq{3tag} simplifies to the product equality
\Eq{*}{
  (q_0-q_i)(y)\cdot(g(x)-g(y))\cdot(q_i(x)-\bar{q}_i(x))=0.
}
The first factor is not zero, because it is the sum of positive terms. Using that $x<\alpha<y_0-\delta_0<y$, the strict monotonicity of $g$ implies that $g(x)<g(y)$, proving that the second factor is also not zero.
Therefore, we must have $q_i(x)=\bar{q}_i(x)$, which shows that $x\in\,]\alpha-\delta_0,\alpha[\,\setminus S_i$. Conversely, if $x\in\,]\alpha-\delta_0,\alpha[\,\setminus S_i$, then $q_i(x)=\bar{q}_i(x)$. In this case \eq{3tag} reduces to the product equality
\Eq{*}{
  q_i(x)\cdot(q_i(x)+(q_0-q_i)(y))\cdot(g(x)-\bar{g}(x))=0.
}
The first two factors are positive, hence we must have $g(x)=\bar{g}(x)$, which proves that $x\in\,]\alpha-\delta_0,\alpha[\,\setminus S$ and completes the proof of the equality \eq{SS}. The maximality of the interval $J$, in view of \eq{SS} implies that
\Eq{SSsup}{
   \sup S\,\cap \,]\alpha-\delta_0,\alpha[\,
   = \sup S_i\,\cap \,]\alpha-\delta_0,\alpha[\,=\alpha \qquad(i\in\{1,\dots,n\}).
}

Let $i\in\{1,\dots,n\}$ be fixed and $y_1,y_2\in\,]y_0-\delta_0,y_0+\delta_0[$ be arbitrary such that $y_1\neq y_2$. Replacing $y$ by $y_1$ and $y_2$ in \eq{3tag}, and then subtracting the two equations so obtained side by side, we get that
\Eq{y1y2}{
    &((q_0-q_i)(y_1)-(q_0-q_i)(y_2))\cdot(q_i(x)g(x)-\bar{q}_i(x)\bar{g}(x))\\
	&\quad+((q_0-q_i)(y_1)g(y_1)-(q_0-q_i)(y_2)g(y_2))\cdot(\bar{q}_i(x)-q_i(x))=0.
}
Let $x_1, x_2\in\,]\alpha-\delta_0,\alpha[\,$ be arbitrary. Substituting $x$ by $x_1$ and then $x_2$ in \eq{y1y2}, we get a homogeneous linear system of two equations of the form
\Eq{*}{
  \xi \cdot(q_i(x_i)g(x_i)-\bar{q}_i(x_i)\bar{g}(x_i))
  +\eta \cdot(\bar{q}_i(x_i)-q_i(x_i))=0 \qquad(i\in\{1,2\}),
}
which is nontrivially solvable with respect to $(\xi,\eta)$, because the equalities 
\Eq{HLE}{
    \xi:=(q_0-q_i)(y_1)-(q_0-q_i)(y_2)=0\quad\mbox{and}\quad 
    \eta:=(q_0-q_i)(y_1)g(y_1)-(q_0-q_i)(y_2)g(y_2)=0
} 
cannot be satisfied simultaneously. Indeed, if $\xi=0$, then $(q_0-q_i)(y_1)=(q_0-q_i)(y_2)>0$. This equality together with $\eta=0$ implies that $g(y_1)=g(y_2)$. The strict monotonicity of $g$ then yields $y_1=y_2$, which contradicts the choice of $y_1$ and $y_2$. Hence the determinant of the system \eq{HLE} must be equal to zero, that is,
\Eq{*}{
	\begin{vmatrix}
		q_i(x_1)g(x_1)-\bar{q}_i(x_1)\bar{g}(x_1) & \bar{q}_i(x_1)-q_i(x_1) \\ q_i(x_2)g(x_2)-\bar{q}_i(x_2)\bar{g}(x_2) & \bar{q}_i(x_2)-q_i(x_2)
	\end{vmatrix}
	=0.
}
If $x_1,x_2\in S\,\cap \,]\alpha-\delta_0,\alpha[\,=S_i\,\cap \,]\alpha-\delta_0,\alpha[\,$ are arbitrary, then $\bar{q}_i(x_1)\neq q_i(x_1)$ and $\bar{q}_i(x_2)\neq q_i(x_2)$, therefore, the above determinantal equality can be rewritten as
\Eq{*}{
  \frac{q_i(x_1)g(x_1)-\bar{q}_i(x_1)\bar{g}(x_1)}{\bar{q}_i(x_1)-q_i(x_1)}
  =\frac{q_i(x_2)g(x_2)-\bar{q}_i(x_2)\bar{g}(x_2)}{\bar{q}_i(x_2)-q_i(x_2)}.
}
Therefore, there exists a real constant $c_i$ such that
\Eq{*}{
	c_i=\frac{q_i(x)g(x)-\bar{q}_i(x)\bar{g}(x)}{\bar{q}_i(x)-q_i(x)}
}
holds for all $x\in S\,\cap \,]\alpha-\delta_0,\alpha[\,$. Solving this equation with respect to $\bar{g}(x)$, we obtain that
\Eq{barg}{
	\bar{g}(x)=\frac{q_i(x)}{\bar{q}_i(x)}(g(x)+c_i)-c_i
} 
is valid for all $x\in S\,\cap \,]\alpha-\delta_0,\alpha[\,$. Subsituting formula \eq{barg} into \eq{3tag}, for all $x\in S\,\cap \,]\alpha-\delta_0,\alpha[\,$ and $y\in\,]y_0-\delta_0,y_0+\delta_0[\,$, we arrive at the equation
\Eq{*}{
    (\bar{q}_i(x)-q_i(x))\cdot\big(q_i(x)(g(x)+c_i)+(q_0-q_i)(y)(g(y)+c_i)\big)=0,
}
which simplifies to the identity
\Eq{*}{
    q_i(x)(g(x)+c_i)=-(q_0-q_i)(y)(g(y)+c_i) 
    \qquad(x\in S\,\cap \,]\alpha-\delta_0,\alpha[\,,\,y\in\,]y_0-\delta_0,y_0+\delta_0[\,).
}
Therefore, there exists a real constant $d_i$ such that
\Eq{*}{
  q_i(x)(g(x)+c_i)=d_i=-(q_0-q_i)(y)(g(y)+c_i) 
    \qquad(x\in S\,\cap \,]\alpha-\delta_0,\alpha[\,,\,y\in\,]y_0-\delta_0,y_0+\delta_0[\,).
}
Using these equalities on the domain indicated, the inequality \eq{gal} implies that
\Eq{alp}{
  g(\alpha)
  \leq\frac{q_i(x)g(x)+(q_0-q_i)(y)g(y)}{q_i(x)+(q_0-q_i)(y)} 
  =\frac{d_i-c_iq_i(x)-d_i-c_i(q_0-q_i)(y)}{q_i(x)+(q_0-q_i)(y)} =-c_i.
}
Therefore, for all $x\in S\,\cap \,]\alpha-\delta_0,\alpha[\,$, we have that $g(x)<g(\alpha)\leq -c_i$, which yields that $d_i<0$ and $q_i(x)=\frac{d_i}{g(x)+c_i}$. This shows that $q_i$ is strictly increasing on $S\,\cap \,]\alpha-\delta_0,\alpha[\,$.
As a consequence of this property, it follows that the equality $q_i(x)(g(x)+c_i)=d_i$
uniquely determines the constants $c_i$ and $d_i$. Indeed, if $q_i(x)(g(x)+c_i')=d_i'$ were also true for all $x\in S\,\cap \,]\alpha-\delta_0,\alpha[\,$ and for some constant $c_i'$ and $d_i'$, then subtracting the two equations side by side, we get $q_i(x)(c_i-c_i')=d_i-d_i'$. If $c_i\neq c_i'$, then this last equality yields that $q_i$ is constant, which contradicts its strict monotonicity. Therefore, $c_i=c_i'$ implying that $d_i=d_i'$ is also valid.

In the final step, instead of a fixed element $y_0\in\,]\alpha,\alpha+r[\,$, we take another arbitrary element $y'\in\,]\alpha,\alpha+r[\,$. Repeating the same argument as above, there exists a positive number $\delta':=\delta(y')$ and real constants $c_i'$, $d_i'$ such that 
\Eq{*}{
  q_i(x)(g(x)+c_i')=d_i'=-(q_0-q_i)(y)(g(y)+c_i') 
    \qquad(x\in S\,\cap \,]\alpha-\delta',\alpha[\,,\,y\in\,]y'-\delta',y'+\delta'[\,).
}
On the set $S\,\cap \,]\alpha-\min(\delta',\delta_0),\alpha[\,$, we have both $q_i(x)(g(x)+c_i)=d_i$ and $q_i(x)(g(x)+c_i')=d_i'$. Due to the uniqueness property, it follows that $c_i'=c_i$ and $d_i'=d_i$. Therefore, 
\Eq{di}{
d_i=-(q_0-q_i)(y)(g(y)+c_i)
}
is valid for all $y\in\,]y'-\delta',y'+\delta'[\,$, in particular, for $y=y'$. The point $y'$ being arbitrary, we can see that \eq{di} holds for all $y\in\,]\alpha,\alpha+r[\,$. Comparing the signs of both sides, we obtain that $g(y)+c_i>0$ for all $y\in\,]\alpha,\alpha+r[\,$. Upon taking the limit $y\to\alpha+0$, it follows that $g(\alpha)+c_i\geq0$. On the other hand, by \eq{alp}, we also have that $g(\alpha)+c_i\leq0$, whence 
$g(\alpha)+c_i=0$ follows. Using that \eq{SSsup} holds, we may also take the limit $x\to\alpha-0$ in the equality
\Eq{*}{
  q_i(x)(g(x)+c_i)=d_i 
    \qquad(x\in S\,\cap \,]\alpha-\delta_0,\alpha[\,),
}
whence we arrive at the equality $d_i=0$, which is the desired contradiction.
\end{proof}

\section{Partial derivatives of Bajraktarević means}

In the next result we determine the partial derivatives of the Bajraktarević means up to third order at diagonal points of $I^n$ under tight regularity assumptions. For instance, as stated below in assertions (1), (2b), (3c), we prove the existence of partial derivatives of the form $\partial_i^m$ only assuming $(m-1)$ times continuous differentiability of $p_i$. 

\Thm{DB}{
	Let $\ell\in\{1,2,3\}$, let $f: I\to\R$ be an $\ell$ times differentiable function on $I$ with a nonvanishing first derivative, and let $p=(p_1,\dots,p_n):I\to\R_{+}^n$. Then we have the following assertions.
	\begin{enumerate}
		\item[(1)] If $\ell=1$, $i\in \{1,\dots,n\}$, and $p_i$ is continuous on $I$, then the first-order partial derivative $\partial_i A_{f,p}$ exists on $\diag(I^n)$ and
		\Eq{*}{
		\partial_i A_{f,p}\circ\Delta_n=\frac{p_i}{p_0}.
		}
		\item[(2a)] If $\ell=2$, $i,j\in \{1,\dots,n\}$ with $i\neq j$, furthermore, $p_i$ and $p_j$ are differentiable on $I$, then the second-order partial derivative $\partial_i\partial_j A_{f,p}$ exists on $\diag(I^n)$ and
		\Eq{*}{
		\partial_i\partial_j A_{f,p}\circ\Delta_n =-\frac{(p_ip_j)'}{p_0^2}-\frac{p_ip_j}{p_0^2}\cdot\frac{f''}{f'}.
		}
		\item[(2b)] If $\ell=2$, $i\in \{1,\dots,n\}$, and $p_i$ is continuously differentiable on $I$, then the second-order partial derivative $\partial_i^2 A_{f,p}$ exists on $\diag(I^n)$ and
		\Eq{*}{
		\partial_i^2 A_{f,p}\circ\Delta_n=2\frac{p_i'(p_0-p_i)}{p_0^2}+\frac{p_i(p_0-p_i)}{p_0^2}\cdot\frac{f''}{f'}.
		}
		\item[(3a)] If $\ell=3$, $i,j,k\in \{1,\dots,n\}$ with $i\neq j\neq k\neq i$, furthermore, $p_i$, $p_j$, and $p_k$ are differentiable on $I$, then the third-order partial derivative $\partial_i\partial_j \partial_k A_{f,p}$ exists on $I^n$ and
		\Eq{*}{\qquad
		\partial_i \partial_j  \partial_k A_{f,p}\circ\Delta_n=2\frac{p_ip_j'p_k'+p_i'p_jp_k'+p_i'p_j'p_k}{p_0^3}+2\frac{(p_ip_jp_k)'}{p_0^3}\cdot\frac{f''}{f'}+\frac{p_ip_jp_k}{p_0^3}\bigg(3\bigg(\frac{f''}{f'}\bigg)^2-\frac{f'''}{f'}\bigg).
                }
        \item[(3b)] If $\ell=3$, $i,j\in \{1,\dots,n\}$ with $i\neq j$, furthermore, $p_i$ is twice differentiable and $p_j$ is differentiable on $I$, then the third-order partial derivative $\partial_i^2\partial_j A_{f,p}$ exists on $I^n$ and
		\Eq{*}{
		\partial_i^2 \partial_j A_{f,p}\circ\Delta_n&=\frac{2p_i'p_j'(2p_i-p_0)+p_j(2(p_i')^2-p_i''p_0)}{p_0^3}+\frac{(2p_i'p_j+p_ip_j')(2p_i-p_0)}{p_0^3}\cdot\frac{f''}{f'}\\&\quad+\frac{p_ip_j}{p_0^3}\bigg((3p_i-p_0)\bigg(\frac{f''}{f'}\bigg)^2-p_i\frac{f'''}{f'}\bigg).
                }
        \item[(3c)] If $\ell=3$, $i\in \{1,\dots,n\}$ and $p_i$ is twice continuously differentiable on $I$, then the third-order partial derivative $\partial_i^3 A_{f,p}$ exists on $\diag(I^n)$ and
		\Eq{*}{
		\partial_i^3 A_{f,p}\circ\Delta_n&=\frac{3(p_0-p_i)\big(p_0p_i''-2(p_i')^2\big)}{p_0^3} +3\frac{p_i'(p_0-2p_i)(p_0-p_i)}{p_0^3}\cdot\frac{f''}{f'}\\&\quad-\frac{p_i(p_0-p_i)}{p_0^3}\bigg(3p_i\bigg(\frac{f''}{f'}\bigg)^2-(p_0+p_i)\frac{f'''}{f'}\bigg).
		}
	\end{enumerate}
}

\begin{proof}
Let $\ell\in\{1,2,3\}$. Assume that $f: I\to\R$ is an $\ell$ times differentiable function on $I$ with a nonvanishing first derivative. We have the following formulas for the derivatives of $f^{-1}$:
\Eq{fs}{
  \big(f^{-1}\big)'=\frac{1}{f'}\circ f^{-1}, \qquad
  \big(f^{-1}\big)''=-\frac{f''}{(f')^3}\circ f^{-1}, \qquad
  \big(f^{-1}\big)'''=\frac{3(f'')^2-f'f'''}{(f')^5}\circ f^{-1}.  
}

In this proof, let $\delta$ denote the extended Kronecker symbol, which, for $i,j,k\in\N$, is defined by:
\Eq{*}{
	\delta_{ij}:=\begin{cases} 1 &\mbox{if } i=j,\\ 
		0 &\mbox{otherwise} \end{cases}
	\qquad\mbox{and}\qquad
	\delta_{ijk}:=\begin{cases} 1 &\mbox{if } i=j=k,\\ 
		0 &\mbox{otherwise}. \end{cases}
} 
Furthermore, in order to make the calculations shorter, we use the notation $R:=R_{f,p}$, where $R_{f,p}$ was defined in \eq{Rfp}. Then $A_{f,p}=f^{-1}\circ R_{f,p}=f^{-1}\circ R$. 

To compute the partial derivatives of $R$, we introduce the notations
\Eq{*}{
P(x_1,\dots,x_n)&:=p_1(x_1)+\dots+p_n(x_n),\\
Q(x_1,\dots,x_n)&:=p_1(x_1)f(x_1)+\dots+p_n(x_n)f(x_n).
}
Then $R\cdot P=Q$ and we have that
\Eq{ids}{
  P\circ\Delta_n=p_0,\qquad
  Q\circ\Delta_n=p_0f, \qquad
  R\circ\Delta_n=f, \qquad\mbox{and}\qquad
  f^{-1}\circ R\circ\Delta_n=\id.
}

To prove the first assertion of the theorem, let $x\in I$ be fixed. Then, using the continuity of $p_i$ and the differentiability of $f$ at $x$, we get
\Eq{1R}{
  (\partial_i R\circ\Delta_n)(x)
  &=\lim_{y\to x}\frac{R(\Delta_n(x)+(y-x)e_i)-R(\Delta_n(x))}{y-x}\\&=\lim_{y\to x} \frac{1}{y-x}
  \bigg(\frac{(p_0(x)-p_i(x))f(x)+p_i(y)f(y)}{p_0(x)-p_i(x)+p_i(y)}-f(x)\bigg)\\
  &=\lim_{y\to x} \frac{p_i(y)}{p_0(x)-p_i(x)+p_i(y)}
  \cdot\frac{f(y)-f(x)}{y-x}
  =\frac{p_if'}{p_0}(x).
}
Therefore, using the standard differentiation rules, the last identity in \eq{ids} and \eq{1R}, we obtain
\Eq{*}{
 \partial_i A_{f,p}\circ\Delta_n
  &=\partial_i \big(f^{-1}\circ R\big)\circ\Delta_n
  =\bigg(\frac{\partial_i R}{f'\circ f^{-1}\circ R}\bigg)\circ \Delta_n
  =\frac{\frac{p_i}{p_0}f'}{f'}=\frac{p_i}{p_0}.
}
This completes the proof of assertion (1).

For the proof of statement (2a), let $i,j\in\{1,\dots,n\}$ with $i\neq j$ be fixed and assume that $p_i$ and $p_j$ are differentiable and $f$ is twice differrentiable on $I$. Then, for all $\alpha,\beta\in\{i,j\}$ with $\alpha\neq\beta$, the partial derivatives $\partial_\alpha$ and $\partial_\alpha\partial_\beta$ of $P$ and $Q$ and hence of $R$ exist at every point in $I^n$. Furthermore, for all $(x_1,\dots,x_n)\in I^n$, we have that 
\Eq{PQ1}{
	\partial_\alpha P(x_1,\dots,x_n)&=p_\alpha'(x_\alpha), &\qquad
	\partial_\alpha Q(x_1,\dots,x_n)&=(p_\alpha f)'(x_\alpha), \\
	\partial_\alpha\partial_\beta P(x_1,\dots,x_n)&=0, &
	\partial_\alpha\partial_\beta Q(x_1,\dots,x_n)&=0.
}
Differentiating the identity $R\cdot P=Q$ with respect to the $j$th and then with respect to the $i$th variable, in view of the equalities in the second line in \eq{PQ1}, it follows that 
\Eq{*}{
	\partial_i\partial_j R\cdot P+\partial_j R\cdot \partial_i P+\partial_i R\cdot \partial_j P=0
}
holds on $I^n$, whence, using \eq{1R} and \eq{PQ1}, we arrive at
\Eq{mix2R}{
	\partial_i\partial_j R\circ\Delta_n
	=\bigg(-\frac{\partial_j R\cdot \partial_i P+\partial_i R\cdot \partial_j P}{P}\bigg)\circ\Delta_n
	=-\frac{p_jf'}{p_0}\cdot p_i'-\frac{p_if'}{p_0}\cdot p_j'
	=-\frac{(p_i p_j)'f'}{p_0^2}.
}
Applying the chain rule, the first two formulas in \eq{fs} and then \eq{ids}, \eq{1R}, \eq{mix2R}, it follows that
\Eq{*}{	
	\partial_i \partial_j A_{f,p}\circ\Delta_n
	&=\Big(\big(\big(f^{-1}\big)''\circ R\big)\cdot \partial_i R\cdot \partial_j R +\big(\big(f^{-1}\big)'\circ R\big)\cdot\partial_i \partial_j R\Big)\circ\Delta_n\\&=
	-\frac{f''}{(f')^3}\cdot\frac{p_if'}{p_0}\cdot\frac{p_jf'}{p_0}+\frac{1}{f'}\cdot\frac{-(p_i p_j)'f'}{p_0^2}=
	-\frac{(p_ip_j)'}{p_0^2}+\frac{p_ip_j}{p_0^2}\cdot\frac{f''}{f'}.
}

To justify assertion (2b), let $x\in I$ be fixed. Let $i\in\{1,\dots,n\}$ and assume that $p_i$ is continuously differentiable and $f$ is twice differentiable on $I$. Then the partial derivative $\partial_i$ of $P$ and $Q$ and hence of $R$ exist at every point in $I^n$. Differentiating the identity $R\cdot P=Q$ with respect to the $i$th variable, we have that $\partial_iR\cdot P+R\cdot\partial_i P=\partial_i Q$, whence 
\Eq{*}{
	\partial_iR=\frac{\partial_i Q-R\cdot\partial_i P}{P}.
}
Using this, we obtain
\Eq{2R}{
	(\partial_i^2 &R\circ\Delta_n)(x)\\
	&=\lim_{y\to x}\frac{\partial_i R(\Delta_n(x)+(y-x)e_i)-\partial_i R(\Delta_n(x))}{y-x}\\
	&=\lim_{y\to x} \frac{1}{y-x}\Bigg(\frac{(p_if)'(y)-\frac{(p_0(x)-p_i(x))f(x)+p_i(y)f(y)}{p_0(x)-p_i(x)+p_i(y)}p_i'(y)}{p_0(x)-p_i(x)+p_i(y)}-\frac{p_i(x)f'(x)}{p_0(x)}\Bigg)\\
	&=\lim_{y\to x}\bigg(\frac{(p_0(x)-p_i(x))p_i'(y)}{(p_0(x)-p_i(x)+p_i(y))^2}\cdot\frac{f(y)-f(x)}{y-x}+\frac{1}{y-x}\bigg(\frac{(p_if')(y)}{p_0(x)-p_i(x)+p_i(y)}-\frac{(p_if')(x)}{p_0(x)}\bigg)\bigg)\\
	&=\frac{(p_0-p_i)p_i'f'}{p_0^2}(x)+\lim_{y\to x}\bigg(\frac1{p_0(x)-p_i(x)+p_i(y)}\cdot\frac{(p_if')(y)-(p_if')(x)}{y-x} \\&\hspace{5cm}-\frac{(p_if')(x)}{(p_0(x)-p_i(x)+p_i(y))p_0(x)} \cdot \frac{p_i(y)-p_i(x)}{y-x}\bigg)\\
	&=\bigg(2\frac{p_i'(p_0-p_i)f'}{p_0^2}+\frac{p_if''}{p_0}\bigg)(x).
}
Applying standard calculus rules, the first two formulas in \eq{fs} and then \eq{ids}, \eq{1R}, \eq{2R}, we conclude
\Eq{*}{
	\partial_i^2 A_{f,p}\circ\Delta_n
	&=\partial_i^2 \big(f^{-1}\circ R\big)\circ\Delta_n
	=\Big(\big(\big(f^{-1}\big)''\circ R\big)\cdot (\partial_i R)^2+\big(\big(f^{-1}\big)'\circ R\big)\cdot \partial_i^2 R\Big)\circ\Delta_n\\&=
	-\frac{f''}{(f')^3}\bigg(\frac{p_if'}{p_0}\bigg)^2+\frac{1}{f'}\bigg(2\frac{p_i'(p_0-p_i)f'}{p_0^2}+\frac{p_if''}{p_0}\bigg)=2\frac{p_i'(p_0-p_i)}{p_0^2}+\frac{p_i(p_0-p_i)}{p_0^2}\cdot\frac{f''}{f'}.
}

To prove assertion (3a), let $i,j,k\in\{1,\dots,n\}$ with $i\neq j\neq k\neq i$ and assume that $p_i$, $p_j$ and $p_k$ are differentiable on $I$. 
Then, for all $\alpha,\beta,\gamma\in\{i,j,k\}$ with $\alpha\neq\beta\neq\gamma\neq\alpha$, the partial derivatives $\partial_\alpha$, $\partial_\alpha\partial_\beta$ and $\partial_\alpha\partial_\beta\partial_\gamma$ of $P$, $Q$ and hence of $R$ exist at every point in $I^n$. Furthermore, for all $(x_1,\dots,x_n)\in I^n$, we have the equalities in \eq{PQ1} and in addition
\Eq{PQ1+}{
	\partial_\alpha\partial_\beta\partial_\gamma P(x_1,\dots,x_n)&=0, \qquad
	\partial_\alpha\partial_\beta\partial_\gamma Q(x_1,\dots,x_n)&=0.
}
Differentiating the identity $R\cdot P=Q$ with respect to the $k$th variable, then with respect to the $j$th variable and then with respect to the $i$th variable, in view of the last two formulas in \eq{PQ1} and \eq{PQ1+}, we get $\partial_i\partial_j\partial_k R\cdot P+\partial_i\partial_j R\cdot \partial_k P+\partial_i\partial_k R\cdot \partial_j P+\partial_j\partial_k R\cdot \partial_i P=0$. Thus, applying the first formula in \eq{mix2R} and \eq{PQ1}, we arrive at
\Eq{mixx3R}{
	\partial_i\partial_j\partial_k R\circ\Delta_n
	&=\bigg(-\frac{\partial_i\partial_j R\cdot \partial_k P
	+\partial_i\partial_k R\cdot \partial_j P
	+\partial_j\partial_k R\cdot \partial_i P}{P}\bigg)\circ\Delta_n\\
	&=\frac{(p_ip_j)'p_k'f'+(p_ip_k)'p_j'f'+(p_jp_k)'p_i'f'}{p_0^3}
	=\frac{2(p_ip_j'p_k'+p_i'p_jp_k'+p_i'p_j'p_k)f'}{p_0^3}.
} 
Hence, using \eq{fs} and then \eq{ids}, \eq{1R}, \eq{mixx3R}, \eq{mix2R}, we obtain
\Eq{*}{
	\partial_i \partial_j \partial_k A_{f,p}\circ\Delta_n
	&=\partial_i \partial_j \partial_k \big(f^{-1}\circ R\big)\circ\Delta_n \\
	&=\Big(\big(\big(f^{-1}\big)'''\circ R\big)\cdot \partial_i R\cdot \partial_j R\cdot \partial_k R
	+\big(\big(f^{-1}\big)'\circ R\big)\cdot \partial_i \partial_j \partial_k R\\
	&\quad+\big(\big(f^{-1}\big)''\circ R\big) \cdot(
	\partial_i R\cdot \partial_j \partial_k R
	+\partial_j R\cdot \partial_i \partial_k R
	+\partial_k R\cdot \partial_i \partial_j R)\Big)\circ\Delta_n\\ 
	&=\frac{3(f'')^2-f'f'''}{(f')^5}\cdot\frac{p_if'}{p_0}\cdot\frac{p_jf'}{p_0}\cdot\frac{p_kf'}{p_0}+\frac{1}{f'}\cdot2\frac{(p_ip_j'p_k'+p_i'p_jp_k'+p_i'p_j'p_k)f'}{p_0^3}\\
	&\quad-\frac{f''}{(f')^3}\bigg(\frac{p_if'}{p_0}\cdot\frac{-(p_jp_k)'f'}{p_0^2}+\frac{p_jf'}{p_0}\cdot\frac{-(p_ip_k)'f'}{p_0^2}+\frac{p_kf'}{p_0}\cdot\frac{-(p_ip_j)'f'}{p_0^2}\bigg)\\
	&=2\frac{p_ip_j'p_k'+p_i'p_jp_k'+p_i'p_j'p_k}{p_0^3}+2\frac{(p_ip_jp_k)'}{p_0^3}\cdot\frac{f''}{f'}+\frac{p_ip_jp_k}{p_0^3}\bigg(3\bigg(\frac{f''}{f'}\bigg)^2-\frac{f'''}{f'}\bigg).
}

To verify assertion (3b), let $i,j\in\{1,\dots,n\}$ with $i\neq j$ and assume that $p_i$ is twice and $p_j$ is once differentiable on $I$. 
Then, for all $\alpha,\beta\in\{i,j\}$ with the assumption $\alpha$ and $\beta$ are not equal to $j$ simultaneously, the partial derivatives $\partial_\alpha$, $\partial_\alpha\partial_\beta$  and $\partial_i^2\partial_j$ of $P$, $Q$ and hence of $R$ exist at every point in $I^n$. Furthermore, for all $(x_1,\dots,x_n)\in I^n$, we have \eq{PQ1}, \eq{PQ1+}, and in addition
\Eq{2P}{
	\partial_i^2 P(x_1,\dots,x_n)&=p_i''(x_i), &\qquad
	\partial_i^2 Q(x_1,\dots,x_n)&=(p_i f)''(x_i), \\
	\partial_i^2\partial_j P(x_1,\dots,x_n)&=0, &
	\partial_i^2\partial_j Q(x_1,\dots,x_n)&=0.
}
Differentiating the equality $R\cdot P=Q$ with respect to the $j$th variable, and then with respect to the $i$th variable twice, using \eq{PQ1} and \eq{2P}, we get 
\Eq{*}{
\partial_i^2\partial_j R\cdot P+\partial_i^2 R\cdot \partial_j P+2\partial_i\partial_j R\cdot \partial_i P+\partial_j R\cdot \partial_i^2 P=0.
}
Thus, applying \eq{2R}, \eq{PQ1}, \eq{mix2R}, \eq{1R}, and \eq{2P}, we arrive at
\Eq{mix3R}{
	\partial_i^2\partial_j R\circ\Delta_n
	&=\bigg(-\frac{\partial_i^2 R\cdot \partial_j P+2\partial_i\partial_j R\cdot \partial_i P+\partial_j R\cdot \partial_i^2 P}{P}\bigg)\circ\Delta_n \\
	&=-\bigg(2\frac{(p_0-p_i)p_i'f'}{p_0^2}+\frac{p_if''}{p_0}\bigg)\cdot \frac{p_j'}{p_0} 
	+2\frac{(p_i p_j)'f'}{p_0^2}\cdot \frac{p_i'}{p_0}
	-\frac{p_jf'}{p_0}\cdot \frac{p_i''}{p_0}\\
	&=\frac{(2p_i'p_j'(2p_i-p_0)+p_j(2(p_i')^2-p_i''p_0))f'-p_0p_ip_j'f''}{p_0^3}.
} 
Therefore, using \eq{fs} and then \eq{ids}, \eq{1R}, \eq{mix3R}, \eq{mix2R}, \eq{2R}, we get 
\Eq{*}{
	\partial_i^2 \partial_j&A_{f,p}\circ\Delta_n
	=\Big(\big(\big(f^{-1}\big)'''\circ R\big)\cdot(\partial_i R)^2\cdot \partial_j R
	+\big(\big(f^{-1}\big)'\circ R\big)\cdot\partial_i^2  \partial_j R\\
	&\qquad\qquad\quad+\big(\big(f^{-1}\big)''\circ R\big)\cdot \big(
	2\partial_i R\cdot \partial_i \partial_j R
	+\partial_j R\cdot \partial_i^2 R\big)\Big)\circ\Delta_n\\ 
	&=\frac{3(f'')^2-f'f'''}{(f')^5}\cdot\frac{(p_if')^2}{p_0^2}\cdot\frac{p_jf'}{p_0}+\frac{1}{f'}\cdot\frac{(2p_i'p_j'(2p_i-p_0)+p_j(2(p_i')^2-p_i''p_0))f'-p_0p_ip_j'f''}{p_0^3}\\
	&\quad-\frac{f''}{(f')^3}\bigg(-2\frac{p_if'}{p_0}\cdot\frac{(p_ip_j)'f'}{p_0^2}+\frac{p_jf'}{p_0}\cdot\frac{p_0(2p_i'f'+p_if'')-(p_ip_i)'f'}{p_0^2}\bigg)\\
	&=\frac{2p_i'p_j'(2p_i-p_0)+p_j(2(p_i')^2-p_i''p_0)}{p_0^3}+\frac{(2p_i'p_j+p_ip_j')(2p_i-p_0)}{p_0^3}\cdot\frac{f''}{f'}\\
	&\quad+\frac{p_ip_j}{p_0^3}\bigg((3p_i-p_0)\bigg(\frac{f''}{f'}\bigg)^2-p_i\frac{f'''}{f'}\bigg),
}
which completes the proof of case (3b).

To prove assertion (3c), let $i\in\{1,\dots,n\}$ and assume that $p_i$ is twice continuously differentiable on $I$. Then the partial derivatives $\partial_i$, $\partial_i^2$ of $P$, $Q$ and hence of $R$ exist at every point in $I^n$. We have that
\Eq{*}{
  \partial_i^2 R
  =\partial_i\bigg(\frac{\partial_i Q-R\cdot\partial_i P}{P}\bigg)
  =\frac{\partial_i^2Q\cdot P-Q\cdot\partial_i^2P-2\partial_i Q\cdot\partial_i P+2R\cdot(\partial_i P)^2}{P^2}.
}
Then, for all $x,y\in I$, we get
\Eq{2Rnd}{
  \partial_i^2 R(\Delta_n(x)+(y-x)e_i)
  &=\frac{(p_if)''(y)((p_0-p_i)(x)+p_i(y))-\big(((p_0-p_i)f)(x)+(p_if)(y)\big)p_i''(y)}{((p_0-p_i)(x)+p_i(y))^2}\\
  &\quad-\frac{2(p_if)'(y) p_i'(y)-2\frac{((p_0-p_i)f)(x)+(p_if)(y)}{(p_0-p_i)(x)+p_i(y)}(p_i')^2(y)}{((p_0-p_i)(x)+p_i(y))^2}\\
  &=\frac{(p_0-p_i)(x)\big(((p_0-p_i)(x)+p_i(y))p_i''(y)-2(p_i')^2(y)\big)}{((p_0-p_i)(x)+p_i(y))^3}(f(y)-f(x))\\
  &\quad+\frac{2(p_0-p_i)(x)(p_i'f')(y)}{((p_0-p_i)(x)+p_i(y))^2}+\frac{(p_if'')(y)}{(p_0-p_i)(x)+p_i(y)}.
}
Therefore, using \eq{2Rnd} and \eq{2R}, the twice continuous differentiability of $p_i$, we obtain that
\Eq{3R}{
	(\partial_i^3 R&\circ\Delta_n)(x)
	=\lim_{y\to x}\frac{\partial_i^2 R(\Delta_n(x)+(y-x)e_i)-\partial_i^2 R(\Delta_n(x))}{y-x}\\
	&=\lim_{y\to x}\frac{1}{y-x}\bigg(\frac{(p_0-p_i)(x)\big(((p_0-p_i)(x)+p_i(y))p_i''(y)-2(p_i')^2(y)\big)}{((p_0-p_i)(x)+p_i(y))^3}(f(y)-f(x))\\
	&\quad+\frac{2(p_0-p_i)(x)(p_i'f')(y)}{((p_0-p_i)(x)+p_i(y))^2}+\frac{(p_if'')(y)}{(p_0-p_i)(x)+p_i(y)}-\bigg(2\frac{p_i'(p_0-p_i)f'}{p_0^2}+\frac{p_if''}{p_0}\bigg)(x)\bigg)\\
	&=\bigg(\frac{3(p_0-p_i)\big(p_0p_i''-2(p_i')^2\big)f'}{p_0^3}+\frac{3(p_0-p_i)p_i'f''}{p_0^2}+\frac{p_if'''}{p_0}\bigg)(x).
}
Hence, applying \eq{fs}, \eq{1R}, \eq{3R}, and \eq{2R}, we conclude
\Eq{*}{
	\partial_i^3 A_{f,p}&\circ\Delta_n
	=\partial_i^3 \big(f^{-1}\circ R\big)\circ\Delta_n\\
	&=((f^{-1})'''\circ R)(\partial_i R)^3+((f^{-1})'\circ R)\cdot\partial_i^3 R+((f^{-1})''\circ R)(3\partial_i R\cdot \partial_i^2 R)\\
	&=\frac{3(f'')^2-f'f'''}{(f')^5}\bigg(\frac{p_if'}{p_0}\bigg)^3+\frac{1}{f'}\bigg(\frac{3(p_0-p_i)\big(p_0p_i''-2(p_i')^2\big)f'}{p_0^3}+\frac{3(p_0-p_i)p_i'f''}{p_0^2}+\frac{p_if'''}{p_0}\bigg)\\
	&\quad-\frac{f''}{(f')^3}\bigg(3\frac{p_if'}{p_0}\bigg(2\frac{p_i'(p_0-p_i)f'}{p_0^2}+\frac{p_if''}{p_0}\bigg)\bigg)\\
	&=\frac{3(p_0-p_i)\big(p_0p_i''-2(p_i')^2\big)}{p_0^3} +3\frac{p_i'(p_0-2p_i)(p_0-p_i)}{p_0^3}\cdot\frac{f''}{f'}\\&\quad-\frac{p_i(p_0-p_i)}{p_0^3}\bigg(3p_i\bigg(\frac{f''}{f'}\bigg)^2-(p_0+p_i)\frac{f'''}{f'}\bigg),
}
which completes the proof of assertion (3c).
\end{proof}

\Lem{1OC}{
	Let $n\geq 2$ and $f,g: I\to\R$ be differentiable functions on $I$ with nonvanishing first derivatives and $i\in\{1,\dots,n\}$. Let $p=(p_1,\dots,p_n): I\to\R_{+}^n$ and $q=(q_1,\dots,q_n): I\to\R_{+}^n$ such that $p_i$ and $q_i$ are continuous on $I$. If $\partial_i A_{f,p}=\partial_i A_{g,q}$ holds on $\diag(I^n)$, then
        \Eq{pq}{
	\frac{q_i}{q_0}=\frac{p_i}{p_0}
        }
holds on $I$.
}

\begin{proof}
In view of \thm{DB}, we have
\Eq{*}{
	\frac{q_i}{q_0}=\partial_i A_{g,q}\circ\Delta_n
	=\partial_i A_{f,p}\circ\Delta_n=\frac{p_i}{p_0}.
}

\end{proof}

\Lem{2OC}{
	Let $n\geq 2$ and $f,g: I\to\R$ be twice differentiable functions on $I$ with nonvanishing first derivatives. Let $p=(p_1,\dots,p_n): I\to\R_{+}^n$ and $q=(q_1,\dots,q_n): I\to\R_{+}^n$ be continuous functions on $I$ and assume that, for all $i\in\{1,\dots,n\}$, \eq{pq} holds on $I$. Let $j,k\in\{1,\dots,n\}$. Then the following two assertions hold.
	\begin{enumerate}[(i)]
		\item Provided that $j\neq k$ and $p_j$, $p_k$, $q_j$, $q_k$ are differentiable functions on $I$, if $\partial_j\partial_k A_{f,p}=\partial_j\partial_k A_{g,q}$ holds on $\diag(I^n)$, then there exists a nonzero constant $\gamma$ such that, for all $i\in\{1,\dots,n\}$, 
		\Eq{ga}{
			q_i^2 g'=\gamma p_i^2 f'
		}
		is valid on $I$.
		\item Provided that $j=k$ and $p_j$, $q_j$ are continuously differentiable functions on $I$, if $\partial_j^2 A_{f,p}=\partial_j^2 A_{g,q}$ holds on $\diag(I^n)$, then there exists a nonzero constant $\gamma$ such that, for all $i\in\{1,\dots,n\}$, \eq{ga} is valid on $I$.
	\end{enumerate} 
}

\begin{proof}
From \lem{1OC} we obtain that $q_i=r_0p_i$ holds for all $i\in\{0,\dots,n\}$. Assume that $j\neq k$. Then, using \thm{DB}, we have that
\Eq{*}{
	\frac{(p_jp_k)'}{p_0^2}+\frac{p_jp_k}{p_0^2}\cdot\frac{f''}{f'}=\partial_j\partial_k A_{f,p}\circ\Delta_n=\partial_j\partial_k A_{g,q}\circ\Delta_n=\frac{(r_0^2p_jp_k)'}{r_0^2p_0^2}+\frac{r_0^2p_jp_k}{r_0^2p_0^2}\cdot\frac{g''}{g'}.
}
Thus, after reduction, we get that
\Eq{r0}{
	\frac{1}{2}\bigg(\frac{f''}{f'}-\frac{g''}{g'}\bigg)=\frac{r_0'}{r_0}
}
is valid on $I$. Hence, there exists $\gamma\in\R\setminus\{0\}$ such that
\Eq{gamma}{
	r_0=\sqrt{\gamma\cdot\frac{f'}{g'}}
}
holds on $I$, whence, using \lem{1OC} again, it follows that, for all $i\in\{1,\dots,n\}$, \eq{gamma} is valid. 

If $j=k$, then with a similar calculation we arrive at the same differential equation for $r_0$.
\end{proof}

For a three times differentiable function $f: I\to\R$ with a nonvanishing first derivative, we introduce its Schwarzian derivative $S(f): I\to\R$ by the following formula:
\Eq{sd}{
	S(f)=\frac{f'''}{f'}-\frac{3}{2}\bigg(\frac{f''}{f'}\bigg)^2.
}
The following lemma plays a basic role in our proofs.

\Lem{sd}{
	Let $f,g: I\to\R$ be three times differentiable functions on $I$ with nonvanishing first derivatives. If $S(f)=S(g)$ is valid on $I$, then there exist $a,b,c,d\in\R$ with $ad\neq bc$ such that $cf+d$ is positive on $I$ and 
\Eq{gf}{
	g=\frac{af+b}{cf+d}
}
holds on I.
}

Our first main result is contained in the following theorem. It completely characterizes the equality of two generalized Bajraktarević means with at least three variables.

\Thm{eq3}{
	Let $n\geq 3$ and $f,g: I\to\R$ be three times differentiable functions on $I$ with nonvanishing first derivatives. Let $p=(p_1,\dots,p_n): I\to\R_{+}^n$ be a continuous function on $I$ and $q=(q_1,\dots,q_n): I\to\R_{+}^n$. Assume that there exist $i,j,k\in\{1,\dots,n\}$ with $i\neq j\neq k\neq i$ such that $p_i,p_j,p_k$ are differentiable functions on $I$. Then the following assertions are equivalent.
	\begin{enumerate}[(i)]
	 \item The $n$-variable generalized Bajraktarević means $A_{f,p}$ and $A_{g,q}$ are identical on $I^n$.
	 \item There is an open subset $U$ of $I^n$ containing $\diag(I^n)$ such that the $n$-variable generalized Bajraktarević means $A_{f,p}$ and $A_{g,q}$ are identical on $U$.
	 \item The function $q$ is continuous, the functions $q_i,q_j,q_k$ are differentiable on $I$, and the equalities 
	 \Eq{*}{
            \partial_\ell A_{f,p}&=\partial_\ell A_{g,q} \qquad(\ell\in\{1,\dots,n-1\}), \\
            \partial_i\partial_j A_{f,p}&=\partial_i\partial_j A_{g,q}, \\
            \partial_i\partial_j\partial_k A_{f,p}&=\partial_i\partial_j\partial_k A_{g,q} 
         }
         hold on $\diag(I^n)$.
	 \item There exist $a,b,c,d\in\R$ with $ad\neq bc$ such that 
	\Eq{*}{
		g=\frac{af+b}{cf+d}\qquad\mbox{and}\qquad q_\ell=(cf+d)p_\ell\qquad (\ell\in\{1,\dots,n\})
	}
	hold on $I$.
	\end{enumerate}
}

\begin{proof}
The implication (i)$\Rightarrow$(ii) is obvious. Applying \lem{reg}, it is also easy to see that assertion (iii) follows from statement (ii). The implication (iv)$\Rightarrow$(i) is a consequence of \thm{BM3}. It remains to prove that assertion (iii) implies statement (iv). 

Without loss of generality, we can assume that $i=1$, $j=2$, and $k=3$. One can easily see that, if $\partial_\ell A_{f,p}=\partial_\ell A_{g,q}$ holds for all $\ell\in\{1,\dots,n-1\}$, then it is also valid for $\ell=n$. Using \lem{1OC}, we have that $q_\ell=r_0p_\ell$ holds for all $\ell\in\{0,\dots,n\}$. Hence, using the equality $q_\ell'=r_0'p_\ell+r_0p_\ell'$, we get that
\Eq{*}{
	&\hspace{-3mm}2\frac{p_1p_2'p_3'+p_1'p_2p_3'+p_1'p_2'p_3}{p_0^3}+2\frac{p_1p_2p_3'+p_1p_2'p_3+p_1'p_2p_3}{p_0^3}\cdot\frac{f''}{f'}+\frac{p_1p_2p_3}{p_0^3}\bigg(3\bigg(\frac{f''}{f'}\bigg)^2-\frac{f'''}{f'}\bigg)\\&=\partial_1\partial_2\partial_3 A_{f,p}\circ\Delta_n=\partial_1\partial_2\partial_3 A_{g,q}\circ\Delta_n
	\\&=2\frac{r_0^3(p_1p_2'p_3'+p_1'p_2p_3'+p_1'p_2'p_3)}{r_0^3p_0^3}+4\frac{r_0'r_0^2(p_1p_2p_3'+p_1p_2'p_3+p_1'p_2p_3)}{r_0^3p_0^3}+6\frac{r_0(r_0')^2p_1p_2p_3}{r_0^3p_0^3}
	\\&\quad+6\frac{r_0^2r_0'p_1p_2p_3}{r_0^3p_0^3}\cdot\frac{g''}{g'}+2\frac{r_0^3(p_1p_2p_3'+p_1p_2'p_2+p_1'p_2p_3)}{r_0^3p_0^3}\cdot\frac{g''}{g'}+\frac{r_0^3p_1p_2p_3}{r_0^3p_0^3}\bigg(3\bigg(\frac{g''}{g'}\bigg)^2-\frac{g'''}{g'}\bigg).
}
Thus, applying \eq{r0} three times, after reduction, it follows that
\Eq{*}{
	\bigg(\frac{f''}{f'}\bigg)^2-\frac{1}{3}\cdot\frac{f'''}{f'}=\frac{1}{2}\bigg(\frac{f''}{f'}-\frac{g''}{g'}\bigg)^2+\bigg(\frac{f''}{f'}-\frac{g''}{g'}\bigg)\frac{g''}{g'}+\bigg(\frac{g''}{g'}\bigg)^2-\frac{1}{3}\cdot\frac{g'''}{g'}	
}
is valid on $I$. Whence we obtain that $S(f)=S(g)$ holds on $I$. Therefore, using \lem{sd}, there exist $a,b,c,d\in\R$ with $ad\neq bc$ such that $cf+d$ is positive and \eq{gf} holds on I. Substituting \eq{gf} into \eq{gamma}, we get that $r_0=\delta (cf+d)$ holds on $I$, where $\delta:=\sqrt{\frac{\gamma}{ad-bc}}>0$. Therefore,
\Eq{*}{
  q_\ell=r_0p_\ell=(\delta cf+\delta d)p_\ell \qquad (\ell\in\{1,\dots,n\}),
}
and 
\Eq{*}{
  g=\frac{af+b}{cf+d}=\frac{\delta af+\delta b}{\delta cf+\delta d},
}
which proves that assertion (iv) holds with the constant vector $(\bar{a},\bar{b},\bar{c},\bar{d}):=\delta \cdot (a,b,c,d)$.
\end{proof}

Our second main theorem has two variants concerning the regularity assumptions and characterizes the equality of generalized two-variable nonsymmetric Bajraktarević means.

\Thm{eq2}{
	Let $f,g: I\to\R$ be three times differentiable functions on $I$ with nonvanishing first derivatives. Let $p=(p_1,p_2): I\to\R_{+}^2$ and $q=(q_1,q_2): I\to\R_{+}^2$ such that $p_1\neq p_2$.  Assume that there exists $i\in\{1,2\}$ such that one of the following regularity conditions is satisfied.
	\begin{enumerate}[(a)]
	 \item $p_i$ is twice continuously differentiable and $p_{3-i}$ is continuous on $I$.
	 \item $p_i$ is twice differentiable and $p_{3-i}$ is once differentiable on $I$.
	\end{enumerate}
        Then the following assertions are pairwise equivalent. 
	\begin{enumerate}[(i)]
	 \item The two-variable generalized Bajraktarević means $A_{f,p}$ and $A_{g,q}$ are identical on $I^2$.
	 \item There is an open subset $U$ of $I^2$ containing $\diag(I^2)$ such that the two-variable generalized Bajraktarević means $A_{f,p}$ and $A_{g,q}$ are identical on $U$.
	 \item[(iv)] There exist $a,b,c,d\in\R$ with $ad\neq bc$ such that 
	\Eq{*}{
		g=\frac{af+b}{cf+d},\qquad q_1=(cf+d)p_1,\qquad \mbox{and}
		\qquad q_2=(cf+d)p_2
	}
	hold on $I$.
	\end{enumerate}
}

\begin{proof}
	The implication (i)$\Rightarrow$(ii) is obvious. The implication (iv)$\Rightarrow$(i) is a consequence of \thm{BM3}. 
	It remains to prove that (ii) implies statement (iv) in both regularity settings. 
	
	Applying \lem{reg}, one can see that we have the following assertions, from statement (ii), under the regularity assumptions (a) and (b) of \thm{eq2}, respectively.
	\begin{enumerate}[(i)]
	 \item[(iii)] The function $q_i$ is twice continuously differentiable, $q_{3-i}$ is continuous on $I$, furthermore
	 \Eq{*}{
            \partial_i A_{f,p}=\partial_i A_{g,q}, \qquad
            \partial_i^2 A_{f,p}=\partial_i^2 A_{g,q}, \qquad\mbox{and}\qquad
            \partial_i^3 A_{f,p}=\partial_i^3 A_{g,q} 
         }
         hold on $\diag(I^2)$.	 \item[(iii)'] The function $q_i$ is twice differentiable, $q_{3-i}$ is once differentiable on $I$, furthermore
	 \Eq{*}{
            \partial_i A_{f,p}=\partial_i A_{g,q}, \qquad
            \partial_i^2 A_{f,p}=\partial_i^2 A_{g,q}, \qquad\mbox{and}\qquad
            \partial_i^2\partial_{3-i} A_{f,p}
            =\partial_i^2\partial_{3-i} A_{g,q} 
         }
         hold on $\diag(I^2)$.
	\end{enumerate}
	
	Without loss of generality, we can assume that $i=1$. Then, using the first equation of (iii) or (iii)' and \lem{1OC}, we have $q_j=r_0p_j$ for all $j\in\{0,1,2\}$. Due to the equality $r_0=q_i/p_i$, it follows that $r_0$ is twice differentiable. Furthermore, by the second equation of assertion (iii) or (iii)' we have that \eq{r0} also holds by the second statement of \lem{2OC}. Observe that, differentiating \eq{r0}, we can obtain that \Eq{r0''}{\frac{r_0''}{r_0}
	=\frac{1}{4}\bigg(2\frac{f'''}{f'}-2\frac{g'''}{g'}-\bigg(\frac{f''}{f'}\bigg)^2+3\bigg(\frac{g''}{g'}\bigg)^2-2\frac{f''}{f'}\cdot\frac{g''}{g'}\bigg).
	}
Under the regularity assumption (a) of \thm{eq2}, the third equality in condition (iii) and formula (3c) of \thm{DB}, yields that
\Eq{3pd}{
	&\hspace{-5mm}-\frac{p_2\big(6(p_1')^2-3p_1''(p_1+p_2)\big)}{p_0^3}-3\frac{p_1'p_2(p_1-p_2)}{p_0^3}\cdot\frac{f''}{f'}-3\frac{p_1^2p_2}{p_0^3}\bigg(\frac{f''}{f'}\bigg)^2+\frac{p_1p_2(2p_1+p_2)}{p_0^3}\cdot\frac{f'''}{f'}
	\\&=\partial_1^3 A_{f,p}\circ\Delta_2
	=\partial_1^3 A_{g,q}\circ\Delta_2\\
	&=-\frac{r_0p_2\big(6r_0^2(p_1')^2+12r_0r_0'p_1p_1'+6(r_0')^2p_1^2-3r_0(p_1+p_2)(r_0p_1''+2r_0'p_1'+r_0''p_1)\big)}{r_0^3p_0^3}\\
	&\quad\hspace{.7mm}-3\frac{r_0^2p_2(r_0p_1'+r_0'p_1)(p_1-p_2)}{r_0^3p_0^3}\cdot\frac{g''}{g'}-3\frac{r_0^3p_1^2p_2}{r_0^3p_0^3}\cdot\bigg(\frac{g''}{g'}\bigg)^2+\frac{r_0^3p_1p_2(2p_1+p_2)}{r_0^3p_0^3}\cdot\frac{g'''}{g'}.	
}
Hence, from \eq{3pd}, using \eq{r0} and \eq{r0''}, it follows that
\Eq{*}{
	&-3\frac{p_1'p_2(p_1-p_2)}{p_0^3}\bigg(\frac{f''}{f'}-\frac{g''}{g'}\bigg)-3\frac{p_1^2p_2}{p_0^3}\bigg(\bigg(\frac{f''}{f'}\bigg)^2-\bigg(\frac{g''}{g'}\bigg)^2\bigg)+\frac{p_1p_2(2p_1+p_2)}{p_0^3}\bigg(\frac{f'''}{f'}-\frac{g'''}{g'}\bigg)
	\\&\quad+3\frac{p_1'p_2(p_1-p_2)}{p_0^3}\bigg(\frac{f''}{f'}-\frac{g''}{g'}\bigg)
	+\frac{3}{2}\cdot\frac{p_1^2p_2}{p_0^3}\bigg(\frac{f''}{f'}-\frac{g''}{g'}\bigg)^2
	\\&\quad-\frac{3}{4}\Cdot\frac{p_1p_2(p_1+p_2)}{p_0^3}\bigg(2\frac{f'''}{f'}-2\frac{g'''}{g'}-\bigg(\frac{f''}{f'}\bigg)^2+3\bigg(\frac{g''}{g'}\bigg)^2-2\frac{f''}{f'}\Cdot\frac{g''}{g'}\bigg) \\
	&\quad+\frac{3}{2}\Cdot\frac{p_1p_2(p_1-p_2)}{p_0^3}\bigg(\frac{f''}{f'}-\frac{g''}{g'}\bigg)\frac{g''}{g'}=0,
}
whence we get
\Eq{*}{
	\frac{1}{2}\cdot\frac{p_1p_2(p_1-p_2)}{p_0^3}\bigg(\frac{f'''}{f'}-\frac{g'''}{g'}\bigg)-\frac{3}{4}\cdot\frac{p_1p_2(p_1-p_2)}{p_0^3}\bigg(\bigg(\frac{f''}{f'}\bigg)^2-\bigg(\frac{g''}{g'}\bigg)^2\bigg)=0,
}
which simplifies to
\Eq{sf-sg}{
	\frac{1}{2}\cdot\frac{p_1p_2(p_1-p_2)}{p_0^3}\big(S(f)-S(g)\big)=0.
}
Using that $p_1\neq p_2$, by continuity, it follows that there exists an open nonempty subinterval $J\subseteq I$ such that $p_1(x)\neq p_2(x)$ holds for $x\in J$. Therefore, the above equation implies that $S(f)=S(g)$ holds on $J$ and hence, by \thm{BM4}, on $I$. Therefore, using \lem{sd}, there exist $a,b,c,d\in\R$ with $ad\neq bc$ such that $cf+d$ is positive and \eq{gf} holds on I. Substituting \eq{gf} into \eq{gamma}, we get that $r_0=\delta (cf+d)$ holds on $I$, where $\delta:=\sqrt{\frac{\gamma}{ad-bc}}>0$. Therefore, with the same argument as at the end of the proof of \thm{eq3}, we can see that assertion (iv) holds with the constant vector $(\bar{a},\bar{b},\bar{c},\bar{d}):=\delta \cdot (a,b,c,d)$.

Under the assumption (b) of \thm{eq2}, the third equality of condition (iii)' and formula (3b) of \thm{DB} imply that
\Eq{3pd+}{
	&\frac{2p_1'p_2'(p_1-p_2)+p_2(2(p_1')^2-p_1''(p_1+p_2)}{p_0^3}+\frac{(2p_1'p_2+p_1p_2')(p_1-p_2)}{p_0^3}\cdot\frac{f''}{f'}+\frac{p_1p_2(2p_1-p_2)}{p_0^3}\bigg(\frac{f''}{f'}\bigg)^2\\
	&-\frac{p_1^2p_2}{p_0^3}\cdot\frac{f'''}{f'}=\partial_1^2\partial_2 A_{f,p}\circ\Delta_n=\partial_1^2\partial_2 A_{g,q}\circ\Delta_n=\frac{2r_0(r_0p_1'+r_0'p_1)(r_0p_2'+r_0'p_2)(p_1-p_2)}{r_0^3p_0^3}\\&+\frac{r_0p_2\big(2\big(r_0^2(p_1')^2+2r_0r_0'p_1p_1'+(r_0')^2p_1^2\big)-r_0(r_0p_1''+2r_0'p_1'+r_0''p_1)(p_1+p_2)\big)}{r_0^3p_0^3}\\
	&+\frac{r_0^2(2p_2(r_0p_1'+r_0'p_1)+p_1(r_0p_2'+r_0'p_2))(p_1-p_2)}{r_0^3p_0^3}\cdot\frac{g''}{g'}+\frac{r_0^3p_1p_2(2p_1-p_2)}{r_0^3p_0^3}\bigg(\frac{g''}{g'}\bigg)^2-\frac{r_0^3p_1^2p_2}{r_0^3p_0^3}\cdot\frac{g'''}{g'}.
}
Hence, from \eq{3pd+}, using \eq{r0} and \eq{r0''}, we arrive at
\Eq{*}{
	&\frac{2p_1p_1'p_2-2p_1'p_2^2+p_1^2p_2'-p_1p_2p_2'}{p_0^3}\bigg(\frac{f''}{f'}-\frac{g''}{g'}\bigg)+\frac{p_1p_2(2p_1-p_2)}{p_0^3}\bigg(\bigg(\frac{f''}{f'}\bigg)^2-\bigg(\frac{g''}{g'}\bigg)^2\bigg)\\
	&\hspace{1.5cm}-\frac{p_1^2p_2}{p_0^3}\bigg(\frac{f'''}{f'}-\frac{g'''}{g'}\bigg)-\frac{2p_1p_1'p_2-2p_1'p_2^2+p_1^2p_2'-p_1p_2p_2'}{p_0^3}\bigg(\frac{f''}{f'}-\frac{g''}{g'}\bigg)\\
	&\hspace{1.5cm}+\frac{1}{4}\Cdot\frac{p_1p_2(p_1+p_2)}{p_0^3}\bigg(2\frac{f'''}{f'}-2\frac{g'''}{g'}-\bigg(\frac{f''}{f'}\bigg)^2+3\bigg(\frac{g''}{g'}\bigg)^2-2\frac{f''}{f'}\Cdot\frac{g''}{g'}\bigg)\\
	&\hspace{1.5cm}-\frac{1}{2}\Cdot\frac{p_1p_2(2p_1-p_2)}{p_0^3}\bigg(\frac{f''}{f'}-\frac{g''}{g'}\bigg)^2
	-\frac{3}{2}\cdot\frac{p_1p_2(p_1-p_2)}{p_0^3}\bigg(\frac{f''}{f'}-\frac{g''}{g'}\bigg)\frac{g''}{g'}=0,
}
whence we have that \eq{sf-sg} holds, thus following a similar train of thought as above, we get assertion (iv).
\end{proof}

\Thm{eq2=}{
	Let $f,g: I\to\R$ be six times differentiable functions on $I$ with nonvanishing first derivatives. Let $p: I\to\R_{+}$ and $q: I\to\R_{+}$ be continuous functions on $I$ and assume that $p$ is three times differentiable on $I$. Then the following assertions are equivalent.
	\begin{enumerate}[(i)]
	 \item The $2$-variable generalized Bajraktarević means $A_{f,(p,p)}$ and $A_{g,(q,q)}$ are identical on $I^2$.
	 \item There is an open subset $U$ of $I^2$ containing $\diag(I^2)$ such that the $2$-variable generalized Bajraktarević means $A_{f,(p,p)}$ and $A_{g,(q,q)}$ are identical on $U$.
	 \item The function $q$ is three times differentiable and the equalities 
	 \Eq{*}{
            \partial_1^j\partial_2^j A_{f,(p,p)}=\partial_1^j\partial_2^j A_{g,(q,q)}
            \qquad (j\in\{1,2,3\})
         }
         hold on $\diag(I^2)$.
	 \item Either there exist $a,b,c,d\in\R$ with $ad\neq bc$ such that 
	\Eq{*}{
		g=\frac{af+b}{cf+d}\qquad\mbox{and}
		\qquad q=(cf+d)p
	}
	hold on $I$ or there exist two polynomials $P$ and $Q$ of at most second degree such that $P$ and $Q$ are positive on $f(I)$ and $g(I)$, respectively, and there exist two constants $\alpha,\beta\in\R$ such that
	\Eq{*}{
	  g=G^{-1}\circ(\alpha F\circ f+\beta), \qquad
	  p=P^{-\frac12}\circ f,\qquad \mbox{and}
	  \qquad q=Q^{-\frac12}\circ g
	}
	hold on $I$, where $F$ and $G$ denote a primitive function of $1/P$ and $1/Q$, respectively.
	\end{enumerate}
}

\begin{proof}
The implication (i)$\Rightarrow$(ii) is obvious. Applying \lem{reg}, it is also easy to see that assertion (iii) follows from statement (ii). The proof of the  implication (iii)$\Rightarrow$(iv) is based on the result of Losonczi \cite{Los99} (who classified the solutions into 1+32 classes) and a recent characterization of the equality of two-variable (symmetric) Bajraktarević means with two-variable quasi-aritmetic means by Páles and Zakaria \cite{PalZak19}. The proof of the implication (iv)$\Rightarrow$(i) is also described in the paper \cite{PalZak19}.
\end{proof}

\def\MR#1{}

%\bibliography{publ,funcequ,temp}
%\bibliographystyle{amsplain}

\end{document}